\documentstyle[psfig,amssymb]{llncs}

\begin{document}

\newcommand{\eps}{\epsilon}
\newcommand{\ord}{{\cal O}}

\title{One-Dimensional Peg Solitaire}

\author{Cristopher Moore \inst{1,2}
\and
David Eppstein \inst{3}}

\institute{
Computer Science Department, University of New Mexico,
Albuquerque NM 87131 {\tt moore@cs.unm.edu}
\and
Santa Fe Institute, 1399 Hyde Park Road, Santa Fe NM 87501
\and
Dept. of Information and Computer Science, University of California,
Irvine, Irvine CA 92697-3425 {\tt eppstein@ics.uci.edu}}

\maketitle

\begin{abstract}
We solve the problem of one-dimensional peg solitaire.  In particular,
we show that the set of configurations that can be reduced to a single
peg forms a regular language, and that a linear-time algorithm exists
for reducing any configuration to the minimum number of pegs.
\end{abstract}

\section{Introduction}

Peg solitaire is a game for one player.  Each move consists of hopping
a peg over another one, which is removed.  The goal is to reduce the
board to a single peg.  The best-known forms of the game take place on
cross-shaped or triangular boards, and it has been marketed as
``Puzzle Pegs'' and ``Hi-Q.''  Discussions and various solutions can
be found in \cite{kraitchik,gardner,hakmem,bcg,beasley}.

In \cite{guy}, Guy proposes one-dimensional peg solitaire as an open
problem in the field of combinatorial games.  Here we show that the
set of solvable configurations forms a regular language, i.e.\ it can
be recognized by a finite-state automaton.  In fact, this was already
shown in 1991 by Plambeck (\cite{chang}, Introduction and Ch.5) and 
appeared as an exercise in a 1974 book of Manna \cite{manna}.  More 
generally, B. Ravikumar showed that the set of solvable configurations 
on rectangular boards of any finite width is regular \cite{ravikumar}, 
although finding an explicit grammar seems to be difficult on boards 
of width greater than 2.

Thus there is little new about this result.  However, it seems not to
have appeared in print, so here it is.

\begin{theorem}  
The set of configurations that can be reduced to a single peg is the
regular language $0^* L 0^*$ where
\begin{eqnarray}
L & = & 1 + 011 + 110 \nonumber \\
& + & 11 (01)^* 
\,\Bigl[ 00 + 00(11)^+ + (11)^+00 + (11)^* 1011 + 1101 (11)^*
\Bigr]\, (10)^* 11 \nonumber \\
& + & 11 (01)^* (11)^* 01 + 10 (11)^* (10)^* 11 .
\label{lang}
\end{eqnarray}
Here $1$ and $0$ indicate a peg and a hole respectively, $w^*$ means
`0 or more repetitions of $w$', and $w^+=ww^*$ means `1 or more
repetitions of $w$.'
\label{thm}
\end{theorem}

\pagebreak

\begin{proof}
To prove the theorem, we follow Leibnitz \cite{bcg} in starting with a
single peg, which we denote
\[ 1 \]
and playing the game in reverse.  The first `unhop' produces
\[ 011 \mbox{ or } 110 \]
and the next
\[ 1101 \mbox{ or } 1011. \]
(As it turns out, $11$ is the only configuration that cannot be
reduced to a single peg without using a hole outside the initial set
of pegs.  Therefore, for all larger configurations we can ignore the
$0$'s on each end.)

We take the second of these as our example.  It has two ends,
$10\ldots$ and $\ldots11$.  The latter can propagate itself
indefinitely by unhopping to the right,
\[ 1010101011. \]
When the former unhops, two things happen; it becomes an end of the
form $11\ldots$ and it leaves behind a space of two adjacent holes,
\[ 110010101011. \]
Furthermore, this is the only way to create a $00$.  We can move the
$00$ to the right by unhopping pegs into it,
\[ 111111110011. \]
However, since this leaves a solid block of $1$'s to its left, we
cannot move the $00$ back to the left.  Any attempt to do so reduces it
to a single hole,
\[ 111111101111. \]
Here we are using the fact that if a peg has another peg to its left,
it can never unhop to its left.  We prove this by induction: assume it
is true for pairs of pegs farther left in the configuration.  Since
adding a peg never helps another peg unhop, we can assume that the two
pegs have nothing but holes to their left.  Unhopping the leftmost peg
then produces $1101$, and the original (rightmost) peg is still
blocked, this time by a peg which itself cannot move for the same
reason.

In fact, there can never be more than one $00$, and there is no need
to create one more than once, since after creating the first one the
only way to create another end of the form $10\ldots$ or $\ldots01$ is
to move the $00$ all the way through to the other side
\[ 111111111101 \]
and another $00$ created on the right end now might as well be the
same one.

We can summarize, and say that any configuration with three or more
pegs that can be reduced to a single peg can be obtained in reverse
from a single peg by going through the following stages, or their
mirror reflection:
\begin{enumerate}
\item We start with $1011$.  By unhopping the rightmost peg, we obtain
$10(10)^*11$.  If we like, we then
\item Unhop the leftmost peg one or more times, creating a pair of
holes and obtaining $11(01)^*00(10)^*11$.  We can then
\item Move the $00$ to the right (say), obtaining
$11(01)^*(11)^*00(10)^*11$.  We can stop here, or
\item Move the $00$ all the way to the right, obtaining
$11(01)^*(11)^*01$, or
\item Fill the pair by unhopping from the left, obtaining
$11(01)^*(11)^*1011(10)^*11$.
\end{enumerate}
Equation~\ref{lang} simply states that the set of configurations is
the union of all of these plus $1$, $011$, and $110$, with as many
additional holes on either side as we like.  Then $0^* L 0^*$ is
regular since it can be described by a regular expression
\cite{hopcroft}, i.e.\ a finite expression using the operators $+$ and
$*$.  \qed \end{proof}

Among other things, Theorem~\ref{thm} allows us to calculate the
number of distinct configurations with $n$ pegs, which is
\[ N(n) = \left\{ \begin{array}{lll}
1 & \hspace{5mm} & n = 1 \\
1 & & n = 2 \\
2 & & n = 3 \\
15 - 7n + n^2 & & n \ge 4, \,n \mbox{ even} \\
16 - 7n + n^2 & & n \ge 5, \,n \mbox{ odd}
\end{array} \right. \]
Here we decline to count $011$ and $110$ as separate configurations,
since many configurations have more than one way to reduce them.

We also have the corollary

\begin{corollary}\label{one-peg-strategy}
There is a linear-time strategy for playing peg solitaire in one
dimension.
\end{corollary}

\begin{proof}
Our proof of Theorem~\ref{thm} is constructive in that it tells us how
to unhop from a single peg to any feasible configuration.  We simply
reverse this series of moves to play the game.  \qed
\end{proof} 

More generally, a configuration that can be reduced to $k$ pegs
must belong to the regular language $(0^* L 0^*)^k$,
since unhopping cannot interleave the pegs coming from different origins
\cite{chang}.  This leads to the following algorithm:

\begin{theorem}
There is a linear-time strategy for reducing any one-dimensional peg
solitaire configuration to the minimum possible number of pegs.
\end{theorem}

\begin{proof}
Suppose we are given a string $c_0c_1c_2\ldots c_{n-1}$ where each
$c_i\in\{0,1\}$.  Let ${\cal A}$ be a nondeterministic finite
automaton (without $\epsilon$-transitions) for $0^*L0^*$, where $A$ is
the set of states in ${\cal A}$, $s$ is the start state, and $T$ is
the set of accepting states.  We then construct a directed acyclic
graph $G$ as follows: Let the vertices of $G$ consist of all pairs
$(a,i)$ where $a \in A$ and $0 \le i \le n$.  Draw an arc from $(a,i)$
to $(b,i+1)$ in $G$ whenever ${\cal A}$ makes a transition from state
$a$ to state $b$ on symbol $c_i$.  Also, draw an arc from $(t,i)$ to
$(s,i)$ for any $t\in T$ and any $0\le i\le n$.  Since $|{\cal A}| =
\ord(1)$, $|G|=\ord(n)$.

Then any path from $(s,0)$ to $(s,n)$ in $G$ consists of $n$ arcs of
the form $(a,i)$ to $(b,i+1)$, together with some number $k$ of arcs
of the form $(t,i)$ to $(s,i)$.  Breaking the path into subpaths by
removing all but the last arc of this second type corresponds to
partitioning the input string into substrings of the form $0^* L 0^*$,
so the length of the shortest path from $(s,0)$ to $(s,n)$ in $G$ is
$n+k$, where $k$ is the minimum number of pegs to which the initial
configuration can be reduced.  Since $G$ is a directed acyclic graph,
we can find shortest paths from $(s,0)$ by scanning the vertices
$(a,i)$ in order by $i$, resolving ties among vertices with equal $i$
by scanning vertices $(t,i)$ (with $t\in T$) earlier than vertex
$(s,i)$.  When we scan a vertex, we compute its distance to $(s,0)$ as
one plus the minimum distance of any predecessor of the vertex.  If
the vertex is $(s,0)$ itself, the distance is zero, and all other
vertices $(a,0)$ have no predecessors and infinite distance.

Thus we can find the optimal strategy for the initial peg solitaire
configuration by forming $G$, computing its shortest path, using the
location of the edges from $(t,i)$ to $(s,i)$ to partition the
configuration into one-peg subconfigurations, and applying
Corollary~\ref{one-peg-strategy} to each subconfiguration.  Since
$|G|=\ord(n)$, this algorithm runs in linear time.
\end{proof}

In contrast to these results, Uehara and Iwata \cite{uehara} showed
that in two or more dimensions peg solitaire is NP-complete.  However,
the complexity of finding the minimum number of pegs to which a $k
\times n$ configuration can be reduced, for bounded $k > 2$, remains
open.

Finally, we note that Ravikumar has proposed an impartial two-player
game, in which players take turns making peg solitaire moves, and
whoever is left without a move loses.  It is tempting to think that
this two-player game might be PSPACE-complete in two or more
dimensions, and polynomial-time solvable in one.

\paragraph{Acknowledgements.}
We thank Aviezri Fraenkel, Michael Lachmann, Molly Rose, B. Sivakumar,
and Spootie the Cat for helpful conversations, and the organizers of
the 2000 MSRI Workshop on Combinatorial Games.

\end{document}